\documentclass[a4paper,11pt]{amsart}
\usepackage{amssymb}

\newtheorem{theorem}{Theorem}
\newtheorem{definition}{Definition}
\newtheorem{example}{Example}
\newtheorem{proposition}[definition]{Proposition}

\newcommand{\Q}{\mathbb{Q}}
\newcommand{\Z}{\mathbb{Z}}
\newcommand{\R}{\mathbb{R}}
\newcommand{\Co}{\mathbb{C}}

\DeclareMathOperator{\Img}{Im}
\DeclareMathOperator{\Real}{Re}
\DeclareMathOperator{\sign}{sign}
\DeclareMathOperator{\Arg}{Arg}
\DeclareMathOperator{\To}{\rightarrow}

\title[Hyperbolic period of the nonholomorphic Eisenstein series]{Hyperbolic period of the nonholomorphic Eisenstein series}
\author{Mariya Vlasenko}
\email{masha.vlasenko@gmail.com}
\begin{document}
\maketitle

\section{Holomorphic lift of the nonholomorphic Eisenstein series from a closed geodesics}

$\Gamma = PSL(2,\Z)$ acts on quadratic forms from the right by 
$$
Q|g(x,y)=Q(ax+by,cx+dy),
$$
so $gQ=Q|g^{-1}$ defines a left action. 

Let $Q=[A,B,C]$ be a quadratic form with real coefficients and $d(Q)=B^2-4AC > 0$. By $C_Q$ 
we denote the circle $A|z|^2+B \Real z+C=0$ in $\Co$. Then $C_{gQ}=g(C_Q)$.

Consider the function $\phi(z,Q)=\frac{Q(z)}{Q'(z)}$. Then

\begin{proposition} For a quadratic form $Q=[A,B,C]$ with real coefficients and $d(Q)=B^2-4AC > 0$

(i) $\phi(Tz,TQ) = \phi(z,Q)$

(ii) $\phi(\cdot,Q) : \Co P^1 \To \Co P^1$ is a double covering, ramified at $\frac {-B \pm i \sqrt{D}}{2 A}$. It maps upper (lower) halfplane to itself

(iii) $\phi(z,Q) = i \Img z$ for $z \in C_Q$ 
\end{proposition}

By (ii) we can define $\phi(z,Q)^s$ for $s \in \Co$ as analytic function of $z$ in either upper or lower halfplane. Indeed, we put $w^s = e^{s \log(w)}$ and $\log(w) = \log|w|+i \Arg(w)$. Here $\Arg(w) \in [-\pi;\pi]$ is the main branch of argument defined on the cut plane $\Co \backslash (-\infty;0]$.

Due to this and (i) we can consider for $s$ with $\Real s > 1$ and $z$ either in upper or in lower halfplane the function
\begin{equation}\label{eisenstein_series_holomorphic_lift}
F_Q(s,z)=\sum_{g \in \langle T \rangle \backslash \Gamma} \phi(gz,gQ)^s.
\end{equation}
From now we restrict to the upper halfplane for brevity.

\begin{proposition}
This series converges absolutely and $F_Q$ is analytic in each variable in domain $\{\Real s > 1\}\times\{\Img z > 0\}$.
\end{proposition}

Recall the nonholomorphic Eisenstein series
$$
E(s,z)=\sum_{m,n \in \Z}^{\;\;\;\;\;\times}\left(\frac {Im z} {|mz+n|^2}\right)^s=\zeta(2 s) \sum_{(m,n)=1}\left(\frac {\Img z} {|mz+n|^2}\right)^s
$$
$$
= \zeta(2 s)\sum_{g \in \langle T \rangle \backslash \Gamma} \Img (gz)^s.
$$
For $z \in C_Q$ we have $gz \in C_{gQ}$ and $\phi(gz,gQ)=i \Img(gz)$ due to (iii). Hence on $C_Q$
$$
E(s,z)=\frac{\zeta(2 s)}{i^s}F_Q(s,z)=\zeta(2 s)e^{-\frac {\pi}2 i s}F_Q(s,z).
$$

\section{Hyperbolic periods}

Suppose we have a function $\Psi(z,Q)$ of a quadratic form $Q$ and a complex argument $z$, analytic in $z$ (always in the upper halfplane) and satisfying $\Psi(gz,gQ)=\Psi(z,Q)$ for $g \in \Gamma$. Examples of such functions are
$$
\Psi(z,Q)=f(z)Q^k(z), \;\;\; f \in M_{2k}(\Gamma)
$$
and
$$
\Psi(z,Q)=F_Q(s,z), \;\;\; Re s > 1.
$$
Indeed, for the modular form $f$ of weight $2 k$ we have 
$$f(gz)Q^k(gz)=f(z)(cz+d)^{2k} Q^k(gz)=f(z)(Q|g)^k(z)=
f(z)(g^{-1}Q)^k(z).
$$
For our function:
$$
F_{gQ}(s,gz)=\sum_{h \in \langle T \rangle \backslash \Gamma} \phi(hgz,hgQ)^s = F_Q(s,z).
$$

Denote by $\Gamma_Q$ the stabilizer of $Q$ in $\Gamma$. For $\gamma \in \Gamma_Q$ and $z_0$ in the upper halfplane we consider 
\begin{equation}\label{hyp_period}
Int(Q,\gamma) = \int_{z_0}^{\gamma z_0} \Psi(z,Q) \frac{dz}{Q(z)}.
\end{equation}

\begin{proposition}
$Int(Q,\gamma)$ doesn't depend on $z_0$, and
$$
Int(Q,\cdot): \Gamma_Q \To \Co
$$
is an additive function, i.e.
$$
Int(Q,\gamma_1\gamma_2)=Int(Q,\gamma_1)+Int(Q,\gamma_2).
$$ 
\end{proposition}
\begin{proof} 
$$
\frac{d}{dz_0} \int_{z_0}^{\gamma z_0} \Psi(z,Q) \frac{dz}{Q(z)} = \frac{\Psi(\gamma z_0,Q)}{Q(\gamma z_0)}\frac{d}{dz_0}(\gamma z_0)-\frac{\Psi(z_0,Q)}{Q(z_0)}
$$
$$
=\frac{\Psi(z_0,\gamma^{-1}Q)}{\gamma^{-1}Q(z_0)}-\frac{\Psi(z_0,Q)}{Q(z_0)}=0
$$
since $\gamma^{-1}Q=Q$. This also implies additivity because 
$$
\int_{z_0}^{\gamma_1\gamma_2z_0}=\int_{z_0}^{\gamma_2 z_0}+\int_{\gamma_2 z_0}^{\gamma_1\gamma_2 z_0}.
$$
\end{proof}

So, the image of $Int(Q,\cdot)$ on the stabilizer $\Gamma_Q$ is an additive subgroup of $\Co$.
\begin{proposition} 
Let $g \in \Gamma$. The image of $Int(Q,\cdot)$ is the same as of $Int(gQ,\cdot)$.
\end{proposition}
\begin{proof}
Note that $\Gamma_{gQ}=g\Gamma_Q g^{-1}$, and
$$
Int(gQ,g \gamma g^{-1})=\int_{z_0}^{g \gamma g^{-1} z_0} \Psi(z,g Q) \frac{dz}{g Q(z)}
$$
(putting $w=g^{-1}z$, $w_0=g^{-1}z_0$)
$$
=\int_{w_0}^{\gamma w_0} \Psi(g w,g Q) \frac{dw}{(g^{-1} g Q)(w)}=\int_{w_0}^{g w_0} \Psi(w,Q) \frac{dw}{Q(w)}=
Int(Q,\gamma).
$$
\end{proof}

So, the image of $Int(Q,\cdot)$ is an additive subgroup of $\Co$ which depends only on the $\Gamma$-orbit of $Q$.

Now we restrict to the case $Q=[A,B,C]$ with $A,B,C$ integers, $(A,B,C)=1$ and nonsquare $D=B^2-4 A C > 0$. Then $\Gamma_Q = \langle \gamma_Q \rangle$ is an infinite cyclic group, generated by $\gamma_Q = \frac v2 I + \frac u2 N_Q$, where $N_Q=\begin{pmatrix}-B&-2C\\2A&B\end{pmatrix}$ and $(v,u)$ is the smallest positive solution of Pell equation $v^2-Du^2=4$. Moreover, $g \gamma_Q g^{-1} = \frac v2 I + \frac u2 g N_Q g^{-1} = \gamma_{gQ}$ since $g N_Q g^{-1} = N_{gQ}$. So, we get a function of a (narrow) class of integer quadratic forms:
$$
[Q]=\Gamma Q \mapsto Int(Q,\gamma_Q).
$$
In next section we see that in the case $\Psi(Q,z)=F_Q(s,z)$ and $D$ is square free, this leads to the zeta function of the correspondent ideal class in $\Q(\sqrt{D})$.

\section{Hecke integral}

For fixed $Q=[A,B,C]$ (integers) with $D=B^2-4AC$ squarefree, we consider the function of $s$ (for $\Real s > 1$)
\begin{equation}\label{zeta}
s \mapsto -\int_{z_0}^{\gamma_Q z_0} F_Q(s,z) \frac{\sqrt D dz}{Q(z)}.
\end{equation}
(That is $- \sqrt D Int(Q,\gamma_Q)$.) We know from the previous section that it depends only on the orbit $[Q]=\Gamma Q$. This orbit corresponds to a narrow class of ideals in $\Q(\sqrt D)$ --- to a class $\beta$ of ideal $ (2 A \Z + (-B + \sqrt D)\Z)y$ where $y$ is any element with $\sign(Ny) = \sign(A)$. So (\ref{zeta})defines a function $\Phi_{\beta}(s)$ for each narrow ideal class $\beta$. 

Since (\ref{zeta}) does't depend on $z_0$, we take $z_0 \in C_Q$. Then $\gamma_Q z_0 \in C_Q$ and we can integrate along $C_Q$ itself. Then we have the following parametrisation of the upper half of $C_Q$ by 
$$
t \mapsto z_t = \frac{x e^t+x' e^{-t}}{e^t+e^{-t}} + i \frac{x - x'}{e^t+e^{-t}}, t \in (-\infty;+\infty)
$$
where $x>x'$ are roots of $Q$, i.e. $Q(x,1)=Q(x',1)=0$, and $-\frac{\sqrt D dz}{Q(z)}=\sign(A) dt$ (Proposition \ref{length_differential} in Appendix). Also on $C_Q$ we have $F_Q(s,z_t)=\frac{e^{\frac{\pi}2 i s}}{\zeta(2 s)}E(s,z_t)$, whence
$$
- \int_{z_0}^{\gamma_Q z_0} F_Q(s,z) \frac{\sqrt D dz}{Q(z)} = \sign(A) \frac{e^{\frac{\pi}2 is}}{\zeta(s)}
\int_{t_0}^{t_0+2 \sign(A) \log \epsilon}E(s,z_t)dt
$$
where $\epsilon = \frac{v + u \sqrt D}2$ with $(v,u)$--- the smallest positive solution of the Pell equation $v^2-Du^2=4$ (Proposition \ref{shift} in Appendix)
$$
= \frac{e^{\frac{\pi}2 is}}{\zeta(2 s)} \int_{t_0}^{t_0+2 \log \epsilon}E(s,z_t)dt.
$$
This integral was already evaluted by Hecke as follows. Suppose for brevity that $D$ is square free, let $\varepsilon$ be a fundamental unit of $\Q(\sqrt D)$~--- the smallest unit $>1$. Then either $\epsilon=\varepsilon^2$ or $\epsilon=\varepsilon$. Now let $b=\Z + x \Z$ be fraction ideal in $\Q(\sqrt D)$, $[b]$ be its class in class group. Then for $c(s)=\int_{-\infty}^{\infty}\frac{dt}{(e^t+e^{-t})^s}$
$$
\zeta([b]^{-1},s)=\sum_{a \in [b]^{-1}}\frac 1{(Na)^s} = \frac{(Nb)^s}2 \sum_{\lambda \in b \slash \varepsilon}\frac1{|\lambda \lambda'|^s} 
$$
$$
= \frac {(x-x')^s}{2 D^{s/2} c(s)} \sum_{\lambda \in b}^{\;\;\;\times}\int_{t_0}^{t_0+2\log{\varepsilon}}\frac{dt}
{(\lambda^2e^t+(\lambda')^2e^{-t})^s}
$$
$$
=\frac 1{2 D^{s/2} c(s)} \int_{t_0}^{t_0+2\log{\varepsilon}} \zeta_{Q_{z_t}}(s) dt
$$
where $Q_z=\frac 1{\Img z}[|z|^2, 2 \Real z, 1]$ is a definite quadratic form of discriminant -4 for each $z$ in upper halfplane. But
$$
\frac 1{Q_z(m,n)}=\frac {\Img z}h{(m z +n)(m \bar z + n)} = \frac{\Img z}{|m z + n|^2}
$$
so $\zeta_{Q_z}(s)=\sum_{m,n}^{\;\;\;\times}\frac 1{Q_z(m,n)^s} = E(s,z)$. Now we have
$$
\int_{t_0}^{t_0+2\log{\varepsilon}} E(z_t,t) dt = 2 c(s) D^{s/2} \zeta([\Z+\Z x]^{-1},s).
$$
Then for expression (\ref{zeta}) we have 
$$
-\int_{z_0}^{\gamma_Q z_0} F_Q(s,z) \frac{\sqrt D dz}{Q(z)} =
2^{f} c(s) e^{\frac {\pi}2 i s} D^{s/2} \frac{\zeta([\Z+\Z x]^{-1},s)}{\zeta(2 s)}
$$
where $\epsilon=\varepsilon^f$ ($f=1$ or 2). L.h.s. is a function of a narrow class, and r.h.s. is a function of correspondent wide class. From this we get that expressions~(\ref{zeta}) are equal for a pair of conjugated narrow classes (in case they differ, i.e. f=1). Then, denoting (\ref{zeta}) by $\Phi_{\beta}(s)$ for a narrow class $\beta$
\begin{theorem}
For any wide ideal class $\alpha$
$$
\sum_{\beta \in \alpha^{-1}} \Phi_{\beta}(s) = 4 c(s)(i \sqrt D)^s \frac{\zeta(\alpha,s)}{\zeta(2 s)}.
$$
where the sum is taken over (one or two) narrow classes $\beta$ inside $\alpha$.
\end{theorem}

\section{Decomposition of hyperbolic period}

Let again $\Psi(z,Q)$ is analytic in $z$ and satisfies $\Psi(gz,gQ)=\Psi(z,Q)$ for $g \in \Gamma$. Let $\beta$ be a narrow class in $\Q(\sqrt D)$ and $Q_0 = Q, Q_1, \dots, Q_r=Q_0$ be a cycle of reduced forms for $\beta$ with
$Q_j=M_j^{-1}Q_{j-1}=Q_{j-1}|M_j$, 
$$
M_j=\begin{pmatrix} m_j&-1\\1&0 \end{pmatrix}, \;\;\; m_j \ge 2.
$$
Then hyperbolic period of $\Psi(z,Q)$ can be decomposed as follows. We put $z_j = M_1 \dots M_j z_0$, so that $z_r = \gamma_Q z_0$. Then
$$
\int_{z_0}^{\gamma_Q z_0}\Psi(z,Q)\frac {d z}{Q(z)} = \sum_{j=1}^r \int_{z_{j-1}}^{z_j}\Psi(z,Q_0)\frac {d z}{Q_0(z)}
$$
(putting $z = M_1 \dots M_{j-1}w$ in $j$-th integral)
$$
= \sum_{j=1}^r \int_{z_0}^{M_j z_0}\Psi(M_1 \dots M_{j-1}w,Q_0)\frac {d w}{Q_{j-1}(w)}
= \sum_{j=1}^r \int_{z_0}^{M_j z_0}\Psi(w,Q_{j-1})\frac {d w}{Q_{j-1}(w)}.
$$
Now since the sum doesn't depend on $z_0$ we can try to tend it to 0. For any $j$ we have $M_j(0)=\infty$. If $\Psi(z,Q)$ is integrable on imaginary halfaxis, we have:
$$
\int_{z_0}^{\gamma_Q z_0}\Psi(z,Q)\frac {d z}{Q(z)} = \sum_{j=0}^{r-1}\int_0^{i \infty}\Psi(z,Q_j)\frac{d z}{Q_j(z)}.
$$

\begin{example}
$\Psi(z,Q) = 1$. Our decomposition leads to the well known formula
$$
\epsilon^2 = \prod_{j=0}^{r-1} \frac{x_j}{x_j'},
$$
where $x_j$ are all numbers with purely periodic continued fraction correspondent to any given narrow class
and $\epsilon$ is a smallest unit $>1$ of positive norm. 
\end{example}
\begin{proof}
Then 
$$
- \sqrt D \int_{z_0}^{\gamma_Q z_0} \frac{d z}{Q (z)} = 2 \log {\epsilon}
$$
and
$$
- \sqrt D \int_0^{\infty}\frac{d z}{Q(z)} = -sign(A)\int_0^{\infty} \frac{(x-x') dz}{(z-x)(z-x')}
$$
$$
= -sign(A) \int_0^{\infty} \left( \frac 1 {z-x} - \frac 1 {z-x'} \right)dz = - sign(A) \log(\frac{z-x}{z-x'}) \Big|_0^{\infty} = sign(A) \log{\frac x{x'}}.
$$
$Q_j$ are reduced, what means that $A_j>0$ and $0 < x_j' < 1 < x_j$. So,
$$
- \sqrt D \int_0^{\infty} \frac{d z}{Q_j(z)} = \log{\frac {x_j} {x_j'}}.
$$
Taking exponent we get the result.
\end{proof}

Take for example $\Q(\sqrt {15})$. Then $\epsilon = 4+\sqrt{15}$ and one can take the cycle
$x_0=(6+\sqrt{15})/3$, $x_1=(6+\sqrt{15})/7$, $x_2 = (8+\sqrt{15})/7$. Then
$$
(4+\sqrt{15})^2=\frac{(6+\sqrt{15})^2(8+\sqrt{15})}{(6-\sqrt{15})^2(8-\sqrt{15})}.
$$

\begin{example}
Let $\Psi(z,Q)=f(z)Q^k(z)$ where $f \in M_{2k}$. Then $\Psi(z,Q)$ is integrable along $i \R_+$ if $f \in S_{2k}$,
and our formula is the usual decomposition of the hyperbolic period via standard periods for the cusp form.
\end{example}

\begin{example}
Let $\Psi(z,Q)=-\sqrt D F_Q(s,z)$. Then not all terms in (\ref{eisenstein_series_holomorphic_lift}) are integrable along $i \R_+$. Terms for $(c,d) \sim (0,1)$, $(2 A, B)$ are not integrable on $\infty$ and terms with $(c,d) \sim (1,0)$, $(B, 2 C)$~--- at 0. The remaining terms can be integrated and their sum is a well defined function of (reduced) quadratic form
$$
\Phi_Q(s)=\sum_{\frac c d \ne \frac B{2A},\frac {2C}B, 0,\infty} \int_0^{i \infty} \left( \frac{\phi(z,Q)}
{(c z+d)(c N_Qz + d)} \right)^s \frac{- \sqrt D dz}{Q(z)}\;\;\; \Real s > 1.
$$
\end{example}
\begin{proof}
Obviosly, $\frac 1{|Q(z)|}$ is integrable (since $A,B,C \ne 0$). So to prove convergence of the above sum we find a summable in power $Re s > 1$ bound for $\left| \frac{\phi(z,Q)} {(c z+d)(c N_Qz + d)} \right|$. We note that $2 \phi(z,Q) = z - N_Qz$, so we need a bound for
$$
\frac 1{c^2} \left| \frac 1{z+\frac d c} - \frac 1{N_Q z+\frac d c}  \right|
\le \frac 1{c^2} \left( \left| \frac 1{z+\frac d c}\right| + \left| \frac 1{N_Q z+\frac d c}  \right| \right),\;\;\;\; z \in i \R_+
$$
Obviolusly $|z + \frac dc| > |\frac dc|$ and $\frac 1{|c d|}$ is summable in any power $>1$. For the second part we see that $N_Q(i \R_+)$ is a semicircle between points $-\frac {2C}B$ and $-\frac B{2A}$. Remember that $A,C>0$and $B<0$ for a reduced form. If $-\frac d c < 0$ then the distance to either of these two points is greater then $|\frac d c|$. If $\frac dc < 0$ is closer to $-\frac {2C}B$, then 
$$
c^2 \left| \frac dc - \frac {2C}B \right| = \frac {|B|}{|c||2C|c|+B|d||}. 
$$
Then for each $c$ we have $\sum^{\;\;\; \times}_{d} \left(\frac {|B|}{|2C|c|+B|d||}\right)^s < 2|B|^s \zeta(s)$ since $d$ is defined by $|2C|c|+B|d||$ in no more than two ways. So, 
$$
\sum_{c,d}^{\;\;\;\times} \left(\frac {|B|}{|c||2C|c|+B|d||}\right)^s < 4 |B|^s \zeta(s)^2.
$$
And analogously for the point $-\frac B{2A}$.
\end{proof}

Finally for the narrow ideal class $\beta$ we have a decomposition of the form
$$
\Phi_{\beta}(s)=\sum_j \Psi_{Q_j}(s) + \underset{z_0 \To \infty}{\lim} \sum_j \int_{z_0}^{M_j z_0}R_{Q_{j-1}}(s,z)\frac{dz}{Q_{j-1}(z)}
$$
where 
$$
R_Q(s,z) = \sum_{\frac c d = \frac B{2A},\frac {2C}B, 0,\infty}\phi(gz,gQ)^s.
$$
We see this limit should exist. But it depends on choices (maybe we should add some other terms there) and we don't know how to define it as a function of the narrow class $\beta$.

\section*{Appendix: the stabilizer of an indefinite quadratic form and Pell equation}

Let $Q=[A,B,C]$ be an indefinite real quadratic form with roots $x' < x$, let $C_Q$ be the semicircle in upper halfplane joining $x'$ and $x$. We associate to $Q$ the matrix $N_Q=\begin{pmatrix}-B&-2C\\2A&B\end{pmatrix}$.
This matrix has the following important properties:\\
(i) $\det N_Q = -D$\\
(ii) $N_Q^2=D$\\
(iii) $g N_Q g^{-1} = N_{g Q}$ for $g \in SL(2,\R)$\\

We will see that for integer $A,B,C$ with $(A,B,C)=1$ and nonsquare $D$ the stabilizer of $Q$ in $SL(2,\Z)$ is
$$
\{ \frac v2 + \frac u2 N_Q | (v,u) \in \Z^2, v^2-Du^2=4 \}.
$$

\begin{proposition} The stabilizer $Stab_Q$ of $Q$ in $PSL(2,\R)$ is
$$
(\pm) M_t = \cosh(t/2)+ \sinh(t/2) \frac 1{\sqrt D} N_Q, \;\;\; t \in \R.
$$
$M_t$ is a hyperbolic shift by $t$ on $C_Q$ (standard hyperbolic metric) in direction from $x'$ to $x$ (from $x$ to $x'$) when $A>0$ ($A < 0$). 
\end{proposition}
\begin{proof}
$g \in Stab_Q$ iff $g x = x$ and $g x' = x'$. We first consider the following degenerate case: $Q=[0,1,0]$. Then $C_Q = i \R_+$, and hyperbolic distance is $d(iy_1,iy_2)=|log \frac{y_1}{y_2}|$. To stabilize 0 and $\infty$ the matrix should be of the form $g=\begin{pmatrix}a&0\\0&1/a\end{pmatrix}$. It is easy to check that
$$
Y_t = \begin{pmatrix}e^{t/2}&0\\0&e^{-t/2}\end{pmatrix} 
= \cosh(t/2) - \sinh(t/2) \begin{pmatrix}-1&0\\0&1\end{pmatrix}
$$
is a hyperbolic shift by $t$ in direction from 0 to $\infty$. 

Now let $Q=[A,B,C]$, $A \ne 0$ with roots $x'<x$. Consider any matrix $g \in SL(2,\R)$ which sends $0,\infty$ to $x', x$ correspondingly. They all are of the form $\begin{pmatrix}ax&bx'\\a&b\end{pmatrix}$ with $ab=\frac 1{x-x'}$. For any such $g$ we have $g[0,1,0]=-\frac 1{x-x'}[1,-(x+x'),xx']$. Then shift by $t$ in direction from $x'$ to $x$ is
$$
X_t = g Y_t g^{-1} = \cosh(t/2) - \sinh(t/2) N_{g[0,1,0]}
$$
(because of (iii))
$$
= \cosh(t/2) + \sinh(t/2) \frac 1{x-x'} \begin{pmatrix}x+x'&-2xx'\\2&-(x+x')\end{pmatrix} 
$$
$$
= \cosh(t/2) + \sinh(t/2) \frac 1{A(x-x')} \begin{pmatrix}-B&-2C\\2A&B\end{pmatrix} 
$$
$$
= \cosh(t/2) + \sign(A) \sinh(t/2) \frac 1{\sqrt D} \begin{pmatrix}-B&-2C\\2A&B\end{pmatrix}
= M_{\sign(A)t}
$$
\end{proof}

From the last proof we can easily get a parametrization of $C_Q$ by hyperbolic length. Indeed, taking 
$g = \frac 1{\sqrt{x-x'}}\begin{pmatrix}x&x'\\1&1\end{pmatrix}$ we have $g(i \R_+)=C_Q$ and
$$
z_t = g(i e^t) = \frac{x e^t+x' e^{-t}}{e^t+e^{-t}} + i \frac{x - x'}{e^t+e^{-t}}
$$
is on the distanse $t$ from $z_0 = g(i) = \frac{x+x'}2 + i \frac{x-x'}2$. In the following proposition
we find a holomorphic lift of the differential of hyperbolic length on $C_Q$.

\begin{proposition}\label{length_differential}
On $C_Q$ we have $- \frac{\sqrt D dz}{Q(z)}=\sign(A) dt$.
\end{proposition}
\begin{proof}
We use parametrization as above and denote $x_t = Re z_t$, $y_t = Im z_t$. Then on $C_Q$
$$ 
\sign(A) \frac{\sqrt D dz_t}{Q(z_t)dt} = \frac{\sign(A) \sqrt D (\dot x_t + i \dot y_t)}{Q(z_t)} 
=\frac{\sign(A) \sqrt D (\dot x_t + i \dot y_t)}{i y_t(2 A x_t + B)-2 A y_t^2}
$$
(since $A |z_t|^2+ B x_t + C = 0$)
$$
= \frac{(x-x')(\dot x_t + i \dot y_t)}{i y_t(2 x_t - x - x')-2 y_t^2}
$$
$$
= \frac{(x-x')((xe^t-x'e^{-t})(e^t+e^{-t})-(xe^t+x'e^{-t})(e^t-e^{-t}) - i (x-x')(e^t-e^{-t}))}
{i(x-x')(2xe^t+2x'e^{-t}-(x+x')(e^t+e^{-t}))-2(x-x')^2}
$$
$$
=\frac{2(x-x')-i(x-x')(e^t-e^{-t})}{i(x-x')(e^t-e^{-t})-2(x-x')}=-1
$$
\end{proof}

Now we switch to the case when $Q$ is an integral primitive (i.e. (A,B,C)=1) quadratic form with nonsquare discriminant $D = B^2-4 AC > 0$. We are interested in $\Gamma_Q = Stab_Q \cap PSL(2,\Z)$.

\begin{proposition}\label{shift}
$M_t \in SL(2,\Z)$ if and only if $t = 2 \log(\frac{v + u  \sqrt D}2)$ where $(v,u)$ is an integral solution of the Pell equation $v^2 - D u^2 = 4$ with $v>0$. 
\end{proposition}
\begin{proof}
Denote $\frac v2 = \cosh(t/2)$, $\frac u2 = \frac{\sinh(t/2)}{\sqrt D}$. Then $v>0$ and 
$$
1=\cosh(t/2)^2-\sinh(t/2)^2=\frac{v^2-D u^2}4.
$$
So, $(v,u)$ is a real solution of Pell equation with $v>0$ and $M_t = \frac v2 + \frac u2 N_Q$.

Now one can see that $M_t \in SL(2,\Z)$ iff $u,v$ are both integers. Indeed, 
$$
M_t = \begin{pmatrix} \frac 12 (v - B u)&-C u \\ A u & \frac 12 (v + Bu)\end{pmatrix}.
$$
Then $Au,Cu, Bu \in \Z$, so $u \in \Z$ since $(A,B,C)=1$. Also $v = tr M_t \in \Z$.
Analogously one can show that $M_t \in SL(2,\Z)$ when $v,u$ are integers.
\end{proof}

So, $\Gamma_Q = \{ \frac v2 + \frac u2 N_Q | (v,u) \in \Z^2, v^2-Du^2=4 \}$. Elements of stabilizer correspond to (totally positive) units in $\Q(\sqrt D)$: if $\epsilon = \frac v2 + \frac u2 D > 0$ then $M(\epsilon) = M_{2 \log(\epsilon)} =  \frac v2 + \frac u2 N_Q$. Due to property (ii) ($N_Q^2=D$) they are also multiplied as units:
$M(\epsilon_1 \epsilon_2)=M(\epsilon_1)M(\epsilon_2)$. So, $\Gamma_Q$ is infinite cyclic group and we fix a generator $\gamma_Q=\frac v2 + \frac u2 N_Q$ where $(v,u)$ is the smallest positive solution of the Pell equation. $\Gamma_Q=\langle \gamma_Q \rangle$. Due to property (iii) for $g \in \Gamma$ we have $g \gamma_Q g^{-1} = \gamma_{g Q}$.
\end{document}